\documentclass[]{article} 
\usepackage[left=2.5cm,centering]{geometry}
\providecommand{\keywords}[1]
{
  \small	
  \textbf{\textit{Keywords -- }} #1
}
\usepackage[english]{babel}
\usepackage{mathtools}
\usepackage{mathrsfs}
\usepackage{cancel}
\usepackage{graphicx}
\usepackage{empheq}
\usepackage{amsfonts,amssymb,amsmath}
\usepackage{enumitem}
\usepackage{authblk}
\usepackage[table]{xcolor}
\usepackage{tikz}
\usepackage{pgfplots}
\usepackage{bm}          
\usepackage{algorithm}
\newtheorem{remark}{Remark}
\newtheorem{proposition}{Proposition}
\newenvironment{proof}{%
\noindent{\bf Proof.\hskip.5em}\ignorespaces}{%
}

\newcommand{\Chi}{\scalebox{1.1}{$\chi$}}
\begin{document}

\title{A Stable Loosely-Coupled Dirichlet-Neumann Scheme for Fluid-Structure Interaction with Large Added Mass
}

\author[1]{Francesca Renzi\thanks{ corresponding author}}
\author[1]{Christian Vergara}

\affil[1]{LaBS, Department of Chemistry, Materials, and Chemical Engineering, Politecnico di Milano, Milano, Italy\\
\texttt{francesca.renzi@polimi.it, christian.vergara@polimi.it}}
\date{}

\maketitle

\begin{abstract}
Solving fluid-structure interaction (FSI) problems when the densities are similar (large added mass), such as in hemodynamics, is challenging since the stability and convergence of the adopted numerical scheme could be compromised.
In particular, while loosely coupled (LC) partitioned approaches are appealing due to their computational efficiency, the stability issues arising in high added mass regimes limit their applicability.

In this work, we present a new strongly-coupled (SC) partitioning strategy for the solution of the FSI problem, from which we derive a stable LC scheme based on Dirichlet and Neumann interface conditions. 
We analyse the convergence of the new SC scheme on a benchmark problem, demonstrating enhanced behaviour over the standard DN method for specific ranges of a parameter $\alpha$, without additional relaxation. 
Building on this, we introduce a new LC scheme by performing a single iteration per time step.
Stability analysis on a benchmark problem proves that the proposed LC scheme is conditionally stable in large added mass regimes, under a constraint on a parameter $\alpha$. 

Numerical experiments in large added mass settings confirm the theoretical results, demonstrating the effectiveness and applicability of the proposed schemes.

\keywords{Coupled problems, absolute stability, Jury's criterion, added mass effect}

\end{abstract}

\section{Introduction}
\label{intro}
The numerical solution of fluid-structure interaction (FSI) problems in the presence of {\sl large added-mass effect}, i.e. when fluid and structure densities are comparable as it happens in hemodynamics \cite{bazilevs2009computational,crosetto2011fluid,bucelli2023stable,fumagalli2023fluid,fumagalli2024novel}, poses significant challenges when a splitting of the monolithic problem is considered for its numerical solution. 

The splitting approach leads to {\sl partitioned} schemes which solve separately the fluid and the structure problems, by suitably exchanging the coupling conditions at the fluid-structure (FS) interface.
Within this class we can distinguish among {\sl strongly-coupled} (SC) schemes \cite{matthies2002partitioned,badia2008fluid,causin2005added,landajuela2017coupling,nobile2001numerical,kuttler2008fixed,badia2009robin,schussnig2022efficient,degroote2010simulation}, where at each time step sub-iterations are introduced until convergence of the interface conditions, and {\sl loosely-coupled} (LC) schemes \cite{gigante2021stability,banks2014analysis,burman2014explicit,bukavc2013fluid,guidoboni2009stable}, where the two sub-problems are solved just once per time step. 
In the presence of added-mass effect, SC schemes generally suffer from convergence issues, whereas LC schemes are prone to stability issues \cite{causin2005added,forster2007artificial}.
More specifically, the possibility of solving the sub-problems just once per time step makes LC schemes very attractive in terms of computational cost.
However, since the coupling conditions are not exactly satisfied, the work exchanged between the two sub-problems is not perfectly balanced and, especially when the added mass is high, this may induce instabilities in the numerical scheme \cite{causin2005added,badia2008fluid}.

In order to overcome this issue, recent studies have proposed conditionally stable LC schemes owing to suitable interface conditions or stabilization terms,  possibly in combination with operator-splitting methods \cite{burman2014explicit,gigante2021stability,fernandez2013incremental,banks2014analysis,bukavc2013fluid,lukavcova2013kinematic}. 
Yet, all the previous LC schemes are based on non-standard interface conditions and/or require a proper a priori tuning of some parameters (e.g. Robin or Nitsche's parameters).

In this context, we present here a new strategy to obtain a stable LC scheme in the presence of high added mass by maintaining standard Dirichlet and Neumann interface coupling conditions and, at the same time, avoiding the tuning of some method's parameters.

As a first step, we introduce an new SC scheme starting from the interpretation of the standard SC Dirichlet-Neumann (DN) scheme as a Richardson method with a block Gauss-Seidel preconditioner and acceleration parameter $\alpha = 1$. 
By considering more general (convergent) choice of $\alpha$ within this framework, we derive a new family of DN-like SC schemes, referred to as {\sl SC-DN-$\alpha$ algorithms}.  
These methods include correction terms that improve the convergence properties of the standard SC-DN method. 
We carry out a 2D convergence analysis of the SC-DN-$\alpha$ scheme for the benchmark problem proposed in \cite{causin2005added}, demonstrating convergence for a specific range of $\alpha$ without the need for relaxation.

Secondly, we propose a new loosely-coupled scheme (denoted by LC-DN-$\alpha$) obtained by considering just one iteration per time step of the SC-DN-$\alpha$ scheme
 and by adding suitable consistency terms.
On the same benchmark problem as above, we demonstrate the conditional absolute stability of the proposed LC-DN-$\alpha$ scheme in the presence of large added-mass, under suitable constraints on $\alpha$. 

All these theoretical findings are supported by numerical experiments in the hemodynamics regime.
Such results confirm the effectiveness and applicability of the proposed SC-DN-$\alpha$ and LC-DN-$\alpha$ schemes for cardiovascular simulations.
\section{Problem settings}\label{sec:problem_settings}
\subsection{Model equations}
\begin{figure}[h]
  \centering
  \includegraphics[width=0.65\linewidth]{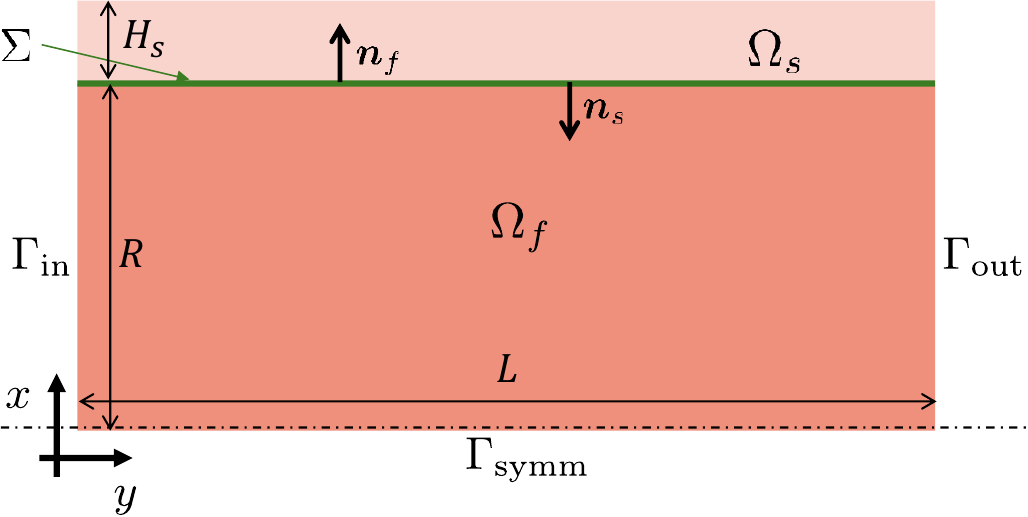}
  \caption{Fluid and structure domains.}
  \label{fig:domain_FSI}
\end{figure}

In this work, we consider the coupling between an incompressible Newtonian fluid and a linear-elastic structure. 
We mathematically described the fluid problem through the incompressible Navier-Stokes (NS) equations formulated in an Arbitrary-Lagrangian-Eulerian (ALE) \cite{donea1982arbitrary} framework and defined on the domain $\Omega^t_f$, with outward normal $\bm{n}_f$. 
The structure problem is described by the linear elasticity equations formulated in a Lagrangian framework, defined on $ \Omega_s^0$ with outward normal $\bm{n}_s$. 
Coupling conditions on the continuity of the fluid and structure velocities and stresses are enforced at the FS interface $\Sigma$. We set $\bm{n} = \bm{n}_f = -\bm{n}_s$. 
A schematic representation of the overall domain is given in Figure \ref{fig:domain_FSI}.

With these settings and given suitable initial conditions at $t = 0$ and boundary conditions on $\partial\Omega_f^t\setminus\Sigma^t$ and $\partial\Omega_s^t\setminus\Sigma^t$, the continuous FSI problem in the strong form reads as follows:

{\sl Find  the fluid velocity $\bm{u}$, the fluid pressure $p$, the structure displacement $\boldsymbol{\eta}$, and the structure interface velocity $\boldsymbol{w}$ satisfying  $\forall t \in (0, T]$, with $T$ the final time, the following problem:}

\begin{align}
\left\{
    \begin{aligned}
        &\rho_f\left[\frac{\delta\bm{u}}{\delta t} + (\bm{u}- \bm{u_{\text{ALE}}})\cdot \nabla\bm{u}\right] - \nabla\cdot\boldsymbol{T}_f(\bm{u}, p) = \boldsymbol{0}
        & \qquad \text{in } \Omega_f^t, \\
        &\nabla \cdot \bm{u} = 0  
        & \qquad \text{in } \Omega_f^t, \\
        &\bm{u} =  \bm{w}
        & \qquad\text{on } \Sigma^t,\\
        &\boldsymbol{T}_f\cdot\bm{n} + \boldsymbol{T}_s\cdot\bm{n} = \boldsymbol{0} 
        & \qquad \text{on } \Sigma^t,\\
        &\rho_s\frac{\partial^2 {\boldsymbol{\eta}}}{\partial t^2} 
        - \nabla\cdot{\boldsymbol{T}}_s ({\boldsymbol{\eta}}) + \beta\bm{\eta} = \boldsymbol{0}
        & \qquad \text{in } {\Omega}_s^0,\\
        &{\bm{w}} = \frac{\partial {\boldsymbol{\eta}}}{\partial t} 
        & \qquad \text{in } {\Sigma}^0,\\
    \end{aligned}
\right.
\label{eq:FSIpb}
\end{align}

\noindent where $\frac{\delta\bm{v}}{\delta t} = \frac{\partial \bm{v}}{\partial t} + (\bm{u}_{\text{ALE}}\cdot \nabla)\bm{v}$ is the ALE time derivative with $\bm{u}_{\text{ALE}}$ the velocity of the fluid domain \cite{donea1982arbitrary,quarteroni2017cardiovascular}, $\beta\neq 0$ surrogates the circumferential elastic forces in the 2D case, whereas $\beta=0$ in the 3D case; 
$\rho_f$ and $\rho_s$ represent the fluid and structure densities.
$\boldsymbol{T}_f$ and $\boldsymbol{T}_s$ are the fluid and the structure Cauchy stress tensors defined as: 
\begin{align*}
    \boldsymbol{T}_f = & -p\boldsymbol{I} + \mu(\nabla\bm{u} + (\nabla\bm{u})^T), \\
    \boldsymbol{T}_s =  & c(\nabla\boldsymbol{\eta} + (\nabla\boldsymbol{\eta})^T) + \lambda (\nabla\cdot\boldsymbol{\eta})\boldsymbol{I},
\end{align*}

\noindent with $\mu$ the dynamic fluid viscosity,  $c$ and $\lambda$  are the Lamé's constants. 
 Notice, for simplicity, zero forcing terms in the previous equations, motivated the main field application we have in mind, i.e., hemodynamics \cite{quarteroni2017cardiovascular}.
Notice also that we do not use two different symbols for the Lagrangian and the Eulerian representation of the unknown. 

\subsection{Time discretization and algebraic setting} \label{subsec:time_discretization}
We adopt the Finite Difference method to discretize in time \eqref{eq:FSIpb}. 
Specifically, we partition the temporal domain $[0,T]$ into the uniform time mesh $t^n = n\Delta t$, with $n = 0, 1, \dots, T/\Delta t$, where $\Delta t$ denotes the time step. 
Given a vector- (or scalar-) valued function $\bm{v}(t)$ (resp. $v(t)$), let $\bm{v}^n$ (resp. $v^n$) denote the approximation of $\bm{v}(t^n)$ (resp. ${v}(t^n)$). 
We linearize the convective term of the NS equations by employing a semi-implicit first-order Euler scheme, and we assume that $\bm{u}_{\text{ALE}}$ and the fluid domain $\Omega_f$ are known at $t^{n+1}$ from extrapolations of previous time steps \cite{nobile2001numerical}. For the structure, we apply the implicit BDF1 scheme.


Since the algorithms presented in this work are motivated by algebraic considerations, we introduce the fully coupled algebraic formulation arising from the space discretization through the Finite Element Method (FEM) of the time-discrete problem \eqref{eq:FSIpb}. 
We start by considering a triangulation of the fluid and structure domains, which generates conforming meshes at the interface. 
We adopt inf-sup stable FE spaces for the fluid.
This leads to the following linear system at each time step: 
\begin{align}
    A\mathbf{X}^{n+1} = \mathbf{b}^{n+1}, 
    \label{eq:FSI-algebric}
\end{align}
where $\mathbf{X}^{n+1} = \left[ \mathbf{U}^{n+1}, \mathbf{P}^{n+1}, \mathbf{U}^{n+1}_{\Sigma},  \mathbf{D}^{n+1}_{\Sigma}, \mathbf{D}^{n+1},\mathbf{W}^{n+1}_{\Sigma},\right]^T$, with $\mathbf{U}^{n+1}$ the vector of the interior fluid velocity, $\mathbf{U}^{n+1}_{\Sigma}$  the vector of the FS interface fluid velocity, $\mathbf{P}^{n+1}$ the vector the fluid pressure, $\mathbf{D}^{n+1}_{\Sigma}$ and $\mathbf{D}^{n+1}$ the vectors of the interface and of the interior structure's displacement, and $\mathbf{W}_{\Sigma}^{n+1}$ the vector of the FS interface structure's velocity. 
The right-hand side vector of \eqref{eq:FSI-algebric},
$\mathbf{b}^{n+1} = \left[ \mathbf{F}^{n+1}_f, \boldsymbol{0}, \boldsymbol{0}, \mathbf{F}^{n+1}_{\Sigma}, \mathbf{F}_s^{n+1}, -\frac{M_{\Sigma}}{\Delta t}{\mathbf{D}^n_{\Sigma}}\right]^T$, collects all the terms arising from the imposition of boundary conditions and time discretizations.
The system matrix is defined as:
\begin{equation*}
A = 
\begin{bmatrix}
	K_{f} & B^T & K_{f,\Sigma} & 0 & 0 & 0\\
	B & 0 & B_{\Sigma} & 0 & 0 & 0 \\
	0 & 0 & M_{\Sigma} & 0 & 0 & \textcolor{red}{- M_{\Sigma}} \\
	K_{\Sigma, f} & (B_{\Sigma})^T & K_{\Sigma\Sigma, f} & K_{\Sigma\Sigma, s} & K_{\Sigma, s} & 0 \\
	0 & 0 & 0 & K_{s, \Sigma} & K_{s} & 0 \\
    0 & 0 & 0 & -\frac{M_{\Sigma}}{\Delta t} &  0 & M_{\Sigma}
\end{bmatrix}, 
\end{equation*}
where ${K}$ are the stiffness matrices, $B$ are associated with the fluid continuity, and $M$ are the mass matrices. Third and fourth block rows are related to the FS interface continuity conditions.
We refer to \cite{badia2008splitting} for details on the derivation of such matrices.

\section{A new strongly coupled partitioned scheme: SC-DN-$\alpha$ method}\label{sec:SC-DN-alpha}

\subsection{Definition of the scheme}\label{subsec:SC-DN-alpha-scheme}

In this section, we present a new staggered algorithm based on sub-iterations between the fluid and the structure sub-problems, where at the FS interface a Dirichlet and a Neumann condition is assigned to fluid and structure sub-problems, respectively. 

We start by rewieving that one iteration of the standard SC-DN scheme is algebraically equivalent to a preconditioned Richardson iteration $P \mathbf{X}^{(k)} = P \mathbf{X}^{(k-1)} + \alpha \left(\mathbf{b} - A \mathbf{X}^{(k-1)}  \right)$, with acceleration parameter $\alpha=1$ and with a block-Gauss-Seidel preconditioner $P$, defined as the matrix $A$ with a zero block, instead of $-M_{\Sigma}$, in position (3, 6), i.e., setting to zero the red term.

 In \cite{causin2005added,forster2007artificial}, it has been proven that this scheme needs a relaxation parameter strictly smaller than 1 to converge when the fluid and solid densities are comparable, which also rapidly approaches zero as the added mass grows larger. 

To overcome this point, we propose here to consider the preconditioned Richardson method applied to \eqref{eq:FSI-algebric}, with $P$ defined as above, associated with
any convergent value of $\alpha$. 
This strategy leads to a new DN-like scheme, named SC-DN-$\alpha$ algorithm and detailed in Algorithm \ref{alg:DN-alpha-SC}. 
We notice that, with respect to the standard SC-DN scheme, \lq\lq correction\rq\rq\, terms appear on the right-hand side.
This will allow us to improve the convergence properties, as detailed in the results reported in the following section.

\begin{algorithm}
\caption{SC-DN-$\alpha$ scheme}
\label{alg:DN-alpha-SC}
{\sl Given suitable initial conditions and boundary conditions on $\partial\Omega_f^t\setminus\Sigma^t$ and $\partial\Omega_s^t\setminus\Sigma^t$ for the fluid velocity and the structure displacement, find, $\forall n = 0,1, \dots, T/ \Delta t$, $\bm{u}^{n+1}, p^{n+1}$, ${\boldsymbol{\eta}}^{n+1}$, and $\bm{w}^{n+1}$ by executing the following steps ($n+1$ understood):}
\vspace{1pt}
\begin{enumerate}[leftmargin=*, align=left]
    \item Set $\boldsymbol{\eta}_{(0)} = \boldsymbol{\eta}^n$, $\bm{w}_{(0)} = \bm{w}^n$, and solve $\forall k = 0, 1, \dots,$ until convergence:
    \vspace{1pt}
    \begin{enumerate}[leftmargin=0.5em, label=\alph*.]
    \item Fluid problem on the unknowns $\bm{u}_{(k+1)} \text{ and } p_{(k+1)}$ equipped with a Dirichlet boundary condition at the interface $\Sigma$: 
    \begin{align}
    \left\{
        \begin{aligned}
           &\rho_f \left(\frac{\bm{u}_{(k+1)}}{\Delta t} + (\bm{u}^n - \bm{u}_{ALE}^n) \cdot \nabla \bm{u}_{(k+1)}\right) - \nabla \cdot \boldsymbol{T}_{f,(k+1)} = \alpha \rho_f \frac{\bm{u}^n}{\Delta t} + \\
           & \qquad + (1- \alpha) \left[\rho_f \left(\frac{\bm{u}_{(k)}}{\Delta t} + (\bm{u}^{n} - \bm{u}_{ALE}^{n}) \cdot \nabla \bm{u}_{(k)}\right) - \nabla \cdot \boldsymbol{T}_{f,(k)}\right] 
            & \text{in}\,\, \Omega_f,\\[6pt]
            &\nabla \cdot \bm{u}_{(k+1)} = (1 - \alpha)\nabla \cdot \bm{u}_{(k)}  
            & \text{in}\,\, \Omega_f,\\[6pt]
            & \bm{u}_{(k+1)} = (1 - \alpha)\bm{u}_{(k)} + \alpha \bm{w}_{(k)}
            & \text{on}\,\, \Sigma;\notag
        \end{aligned}
    \right.
\end{align}

    \item Structure problem on the unknown $\boldsymbol{\eta}_{(k+1)}$ equipped with a Neumann boundary condition at the interface $\Sigma$: 
    \begin{align}
    \left\{
    \begin{aligned}
            &\rho_s \frac{{\boldsymbol{\eta}}_{(k+1)}}{\Delta t^2} - \nabla \cdot {\boldsymbol{T}}_{s,(k+1)} + \beta\bm{\eta}_{(k+1)}= 
            \alpha \rho_s\frac{2{\boldsymbol{\eta}}^n - {\boldsymbol{\eta}}^{n-1}}{\Delta t^2} + \notag \\
            & \qquad\qquad\qquad+ (1-\alpha) \left( \rho_s \frac{{\boldsymbol{\eta}}_{(k)}}{\Delta t^2} - \nabla \cdot {\boldsymbol{T}}_{s,(k)}\right)
            & \text{in}\ {\Omega}_s^0,\\[6pt]
             &\boldsymbol{T}_{s,(k+1)} \cdot \bm{n} = \boldsymbol{T}_{f,(k+1)} \cdot \bm{n}  + (1-\alpha) \left(\boldsymbol{{T}}_{s,(k)} \cdot \bm{n} - \boldsymbol{T}_{f,(k)} \cdot \bm{n} \right)
            & \text{on}\ \Sigma^0; \notag \\[6pt]
    \end{aligned}
    \right.
    \end{align}
    \item Update $\displaystyle{\bm{w}_{(k+1)} = \frac{\bm{\eta}_{(k+1)}-\bm{\eta}^n}{\Delta t}\bigg|_{\Sigma^0}}$;
    \item Check the convergence criterion based on the Euclidean norm of the increment between $k$ and $k+1$ of four unknowns; 
    \end{enumerate}

    \item At convergence, set $\bm{u}^{n+1} = \bm{u}_{(k+1)}, p^{n+1} = p_{(k+1)}$, $\boldsymbol{\eta}^{n+1} = \boldsymbol{\eta}_{(k+1)}$, and $\bm{w}^{n+1} = \bm{w}_{(k+1)}$.
    
\end{enumerate}
\end{algorithm}

\begin{remark}
Notice that Algorithm \ref{alg:DN-alpha-SC} with $\alpha=1$ corresponds to the standard SC-DN scheme.
\end{remark}


\subsection{Convergence of the SC-DN-$\alpha$ algorithm}\label{subsec:convergence}
We present in what follows a theoretical result on the convergence of the SC-DN-$\alpha$ scheme presented in Algorithm \ref{alg:DN-alpha-SC} for a benchmark problem which was introduced in \cite{causin2005added} and later exploited also in \cite{badia2008fluid,gigante2021stability}.

\begin{figure}[h]
  \centering
  \includegraphics[width=0.65\linewidth]{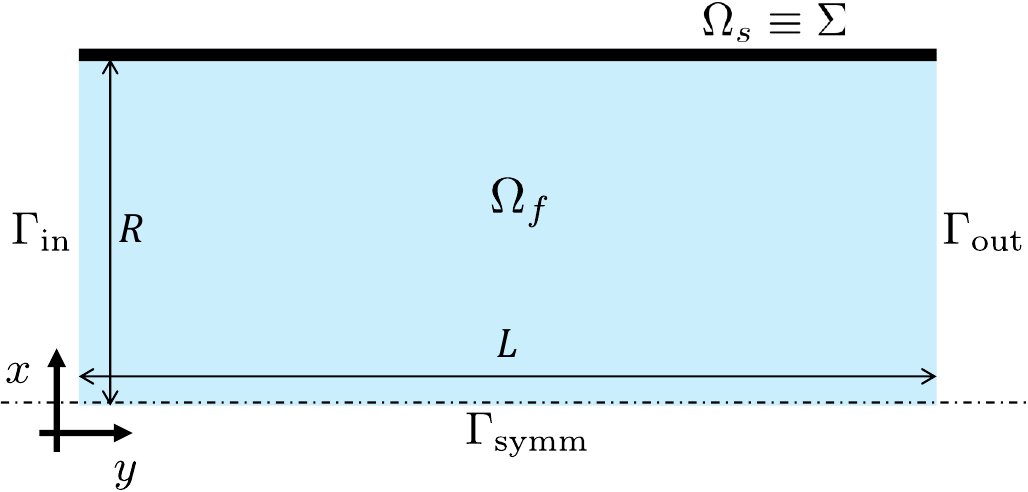}
  \caption{Fluid and structure domains for the benchmark problem.}
  \label{fig:domain}
\end{figure}

In this particular benchmark problem, the fluid domain consists of a fixed 2D rectangular grid $\Omega_f$ with dimension $R\times L$ and boundary $\partial\Omega_f = \Gamma_{\text{in}} \cup \Gamma_{\text{out}} \cup \Gamma_{\text{symm}} \cup \Sigma$, where $\Sigma$ represents the fixed FS interface. 
The structure domain ${\Omega}_s$ is a 1D membrane and thus coincides with  $\Sigma$ (see Figure \ref{fig:domain}).
Moreover, the fluid is modelled as a linear, incompressible, and inviscid fluid, while the structure is modelled through a 1D independent rings model \cite{ma1992numerical,quarteroni2000computational}. 
Regarding the structure displacement, it is assumed to happen only in the direction normal to the interface $\Sigma$.

Unlike in \cite{causin2005added}, we approximate this benchmark problem in time by adopting the backwards Euler for the fluid and the implicit BDF1 scheme for the structure: 

\begin{align}\label{eq:toy_discrete}
    \left\{ 
    \begin{aligned}
        &\rho_f\frac{\bm{u}^{n+1} -\bm{u}^n}{\Delta t} + \nabla p^{n+1} = \bm{0} & \quad \text{in } \Omega_f, \\
        &\nabla\cdot\bm{u}^{n+1} = 0 & \quad \text{in } \Omega_f, \\
        &\bm{u}^{n+1}\cdot\bm{n} = w^{n+1} & \quad \text{on } \Sigma,\\
        &\rho_sh_s\frac{\eta^{n+1} - 2\eta^{n} + \eta^{n-1}}{\Delta t^2} + \beta\eta^{n+1} = p^{n+1} & \quad \text{in } \Omega_s,\\
        & w^{n+1} = \frac{\eta^{n+1}- \eta^{n}}{\Delta t} & \quad \text{in } \Omega_s,
    \end{aligned}
        \right.
\end{align}
with suitable boundary conditions, see \cite{causin2005added}, where $h_s$ is the structure thickness parameter.
Notice that, due to the membrane nature of the structure, the latter equation also serves as the dynamic (third Newton law) condition. 

Accordingly, Algorithm \ref{alg:DN-alpha-SC} applied to the resulting discretized-in-time benchmark problem reads as follows:

{\sl Find $\bm{u}^{n+1}, p^{n+1}$, ${\eta}^{n+1}$, and $w^{n+1}$ by executing the following steps:}
\begin{enumerate}
    \item Set $\eta_{(0)} = {\eta}^n$, $w_{(0)} = w^n$ and solve $\forall k = 0,1,\dots$ until convergence:
    \begin{enumerate}
        \item Fluid problem equipped with a Dirichlet boundary condition at the interface $\Sigma$: 
        \begin{align}
        \left\{
            \begin{aligned}\label{eq:fluid-step-modelSC}
                &\frac{\rho_f}{\Delta t}\bm{u}_{(k+1)} + \nabla p_{(k+1)} = (1 - \alpha)\left(\frac{\rho_f}{\Delta t}\bm{u}_{(k)} + \nabla p_{(k)}\right) 
                & \quad \text{in } \Omega_f, \\
                &\nabla\cdot\bm{u}_{(k+1)} = (1 - \alpha)\nabla\cdot \bm{u}_{(k)}
                & \quad \text{in } \Omega_f, \\
                &\bm{u}_{(k+1)}\cdot\bm{n} = (1-\alpha)\bm{u}_{(k)}\cdot\bm{n} + \alpha w_{(k)}
                & \quad \text{on } \Sigma,
            \end{aligned}
        \right.
        \end{align}
        completed with suitable conditions on the boundary $\partial\Omega_f\backslash\Sigma$; 
        \item  Structure problem: 
        \begin{align}
        \begin{aligned}\label{eq:structure-step_modelSC}
          &  \frac{\rho_s h_s}{\Delta t^2}\eta_{(k+1)} + \beta\eta_{(k+1)} - p_{(k+1)} 
            = \alpha\frac{\rho_s h_s}{\Delta t^2}(2\eta^n - \eta^{n-1}) \\ 
            & \qquad\qquad\qquad + (1 - \alpha)\left(\frac{\rho_s h_s}{\Delta t ^2}\eta_{(k)} + \beta\eta_{(k)} - p_{(k)}\right) 
        \end{aligned}
        \quad \text{in } \Omega_s;\quad\quad
        \end{align}
        \item Update $\displaystyle{w_{(k+1)} = \frac{\eta_{(k+1)} - \eta^n}{\Delta t }}$;
        \item Check the convergence criterion. 
    \end{enumerate}

    \item At convergence, set $\bm{u}^{n+1} = \bm{u}_{(k+1)}$,  ${p}^{n+1} = {p}_{(k+1)}$, $\eta^{n+1} = \eta_{(k+1)}$, and $w^{n+1} = w_{(k+1)}$.
\end{enumerate}

\begin{proposition} The SC-DN-$\alpha$ \eqref{eq:fluid-step-modelSC}-\eqref{eq:structure-step_modelSC} algorithm converges to the solution of \eqref{eq:toy_discrete} if and only if: 
\begin{align}
    0 < \alpha < \frac{2 (\rho_s h_s + \Delta t^2 \beta)}{\rho_s h_s + \Delta t^2\beta + \rho_f\mu_i}.
    \label{eq:prep-thesis}
\end{align}
\end{proposition}

\begin{proof}
Let us consider the momentum equation of \eqref{eq:fluid-step-modelSC} and write it as: 
\begin{align}
    \nabla p_{(k+1)} - (1 - \alpha)\nabla p_{(k)}
    = - \frac{\rho_f}{\Delta t }\left[\bm{u}_{(k+1)} - (1 -\alpha)\bm{u}_{(k)}\right] + \alpha\frac{\rho_f}{\Delta t}\bm{u}^n \quad \text{in } \Omega_f. \qquad
    \label{eq:momentumpressure}
\end{align}

\noindent Multiply this equation by $\bm{n}$ and let $u$ denote the product $ \bm{u} \cdot \bm{n}$: 
\begin{align}
     \frac{\partial p_{(k+1)}}{\partial \bm{n}} - (1 - \alpha) \frac{\partial p_{(k)}}{\partial \bm{n}} = -\frac{\rho_f}{\Delta t}\left(u_{(k+1)} - (1 - \alpha)u_{(k)}\right) + \alpha\frac{\rho_f}{\Delta t }u^n \quad \text{on } \Sigma.
     \qquad
     \label{eq:projectedMomentum}
\end{align}

\noindent Let now write the continuity equation of \eqref{eq:fluid-step-modelSC} as: 
\begin{align*}
    u_{(k+1)} - (1 - \alpha)u_{(k)} = \alpha\frac{\eta_{(k)} - \eta^n}{\Delta t} \quad \text{on}\,\, \Sigma,
\end{align*} 

\noindent and plug it in \eqref{eq:projectedMomentum}: 
\begin{align}
    \frac{\partial p_{(k+1)}}{\partial \bm{n}} - (1 - \alpha)\frac{\partial p_{(k)}}{\partial \bm{n}} = - \alpha\frac{\rho_f}{\Delta t^2}(\eta_{(k)} - \eta^n) + \alpha\frac{\rho_f}{\Delta t}u^n \quad \text{on}\,\, \Sigma.
     \label{eq:substitution}
\end{align}

\noindent Notice that we can write the continuity equation \eqref{eq:fluid-step-modelSC} also as:
\begin{align}
    \nabla \cdot (\bm{u}_{(k+1)} - (1 - \alpha)\bm{u}_{(k)}) = 0 \quad \text{in} \,\, \Omega_f.
    \label{eq:nulldivergence}
\end{align}

\noindent By applying $\nabla \cdot$ to \eqref{eq:momentumpressure} and considering \eqref{eq:nulldivergence} we obtain: 
\begin{align}
    \Delta (p_{(k+1)} &- (1-\alpha)p_{(k)}) = - \frac{\rho_f}{\Delta t}\left(\bm{u}_{(k+1)} - (1 - \alpha)\bm{u}_{(k)}\right)  + \alpha\frac{\rho_f}{\Delta t }\bm{u}^n = 0 \quad \text{in} \,\, \Omega_f.\qquad
    \label{eq:pressure-poisson}
\end{align}

$\\$ Introducing the added mass operator $\mathcal{M}_A : H^{-1/2}(\Sigma) \to H^{1/2}(\Sigma)$ associated with problem \eqref{eq:substitution}-\eqref{eq:pressure-poisson}, which is compact, self-adjoint and positive on $L^2(\Sigma)$ \cite{causin2005added}, we obtain:
\begin{equation*}\label{eq:expressionp}
    p_{(k+1)} - (1-\alpha) p_{(k)} = \mathcal{M}_A \left( - \alpha \frac{\rho_f}{\Delta t^2}(\eta_{(k)} - \eta^{n}) + \alpha\frac{\rho_f}{\Delta t}u^n \right) \quad \text{on}\ \Sigma.
\end{equation*} 

\noindent Inserting the latter in \eqref{eq:structure-step_modelSC}, we obtain:
\begin{align*}
    &\frac{\rho_s h_s}{\Delta t^2} \eta_{(k+1)} + \beta\eta_{(k+1)} - (1 - \alpha)\left(\frac{\rho_s h_s}{\Delta t^2}\eta_{(k)} + \beta\eta_{k}\right) +\\
    &\qquad\qquad +\alpha \frac{\rho_s h_s}{\Delta t^2}(2\eta^{n} - \eta^{n-1}) - \mathcal{M}_A \left( - \alpha \frac{\rho_f}{\Delta t^2}(\eta_{(k)} - \eta^{n}) + \alpha\frac{\rho_f}{\Delta t}u^n \right) = 0.
\end{align*}

\noindent Let $\mu_i = \frac{L}{i\pi\tanh(\frac{i\pi R}{L})}$ be the eigenvalues of $\mathcal{M}_A$ \cite{causin2005added,gigante2021stability}. Due to the compactness of $\mathcal{M}_A$, there exists an orthonormal basis of $L^2(\Sigma)$, made of eigenvectors ${z_i}$ of $\mu_i$, such that we can expand the solution $\eta$ on this basis: $\eta = \sum_i \eta_i z_i$. Each Fourier coefficient $\eta_i$ of the solution must satisfy: 

\begin{align}
    \left(\frac{\rho_s h_s}{\Delta t^2} + \beta\right)\eta_{i,(k+1)}
    &+ \left(\mu_i\alpha \frac{\rho_f}{\Delta t^2} - (1-\alpha)\left[\frac{\rho_s h_s}{\Delta t^2} + \beta\right]\right)\eta_{i,(k)} \notag\\
    &+ g(\eta^n,\eta^{n-1},u^n) = 0,
\label{eq:fixed-point}
\end{align}

\noindent where $g$ collects quantities from previous time steps. 
The necessary and sufficient condition for the convergence of this fixed-point iterative method is: 
\begin{align*}
    \left|\frac{(1-\alpha)(\rho_s h_s + \beta\Delta t^2) - \rho_f\alpha\mu_i}{\rho_s h_s + \beta \Delta t^2}\right| < 1.
\end{align*}
A straightforward calculation leads to \eqref{eq:prep-thesis}. $\qquad \square$
\end{proof}

\begin{remark}
   We point out that the proposed SC-DN-$\alpha$ scheme does not require any relaxation parameter $\omega$ to converge. 
   Moreover, the choice of $\alpha$, unlike $\omega$ in the standard SC-DN scheme, could avoid manual tuning and could exploit the machinery from Richardson's theory.
\end{remark}
\subsection{Numerical results}\label{subsec:convergence-results}
In this section, we present some 2D numerical results obtained by testing the SC-DN-$\alpha$ algorithm (Algorithm \ref{alg:DN-alpha-SC}), with realistic equations, with the aim of validating the theoretical findings reported in the previous section. 
In all the numerical simulations, referring to Figure 1, the structure and the fluid domains are rectangular, with dimensions $L = 6cm, \, R=0.5cm,\,\text{and}\, H_s=0.1cm$, representing a longitudinal section of a cylindrical blood vessel.
We consider the FSI problem \eqref{eq:FSIpb} with $\mu = 0.035poise, \, \rho_f = 1.0g/cm^3, \, \rho_s = 1.1 g/cm^3, \, c = 5\times10^5 dyne/cm^2, \, \lambda = 7.5\times 10^6 dyne/cm^2, \, \beta=5.7\times10^6 dyne/cm^4$, unless otherwise specified. 

Regarding the numerical settings, we adopt the time discretization schemes detailed in Section \ref{subsec:time_discretization}, with time step $\Delta t = 10^{-3}s$ (unless otherwise specified), and in space we use $\mathbb{P}_1-iso\mathbb{P}_2/\mathbb{P}_1$ Finite Elements for the fluid and $\mathbb{P}_1$ for structure variables with a discretization step $h = 0.2cm$.
In all the tests, the tolerance in the normalized stopping criterion is set to $10^{-4}$.
The simulations are carried out using a 2D Finite Element code written in Matlab at MOX, Dipartimento di Matematica of Politecnico di Milano, at CMCS-EPFL-Lausanne, and at LaBS, Dipartimento di Chimica, Materiali e Ingegneria Chimica of Politecnico di Milano. 

\begin{figure}[ht]
  \centering
  \includegraphics[width=1\linewidth]{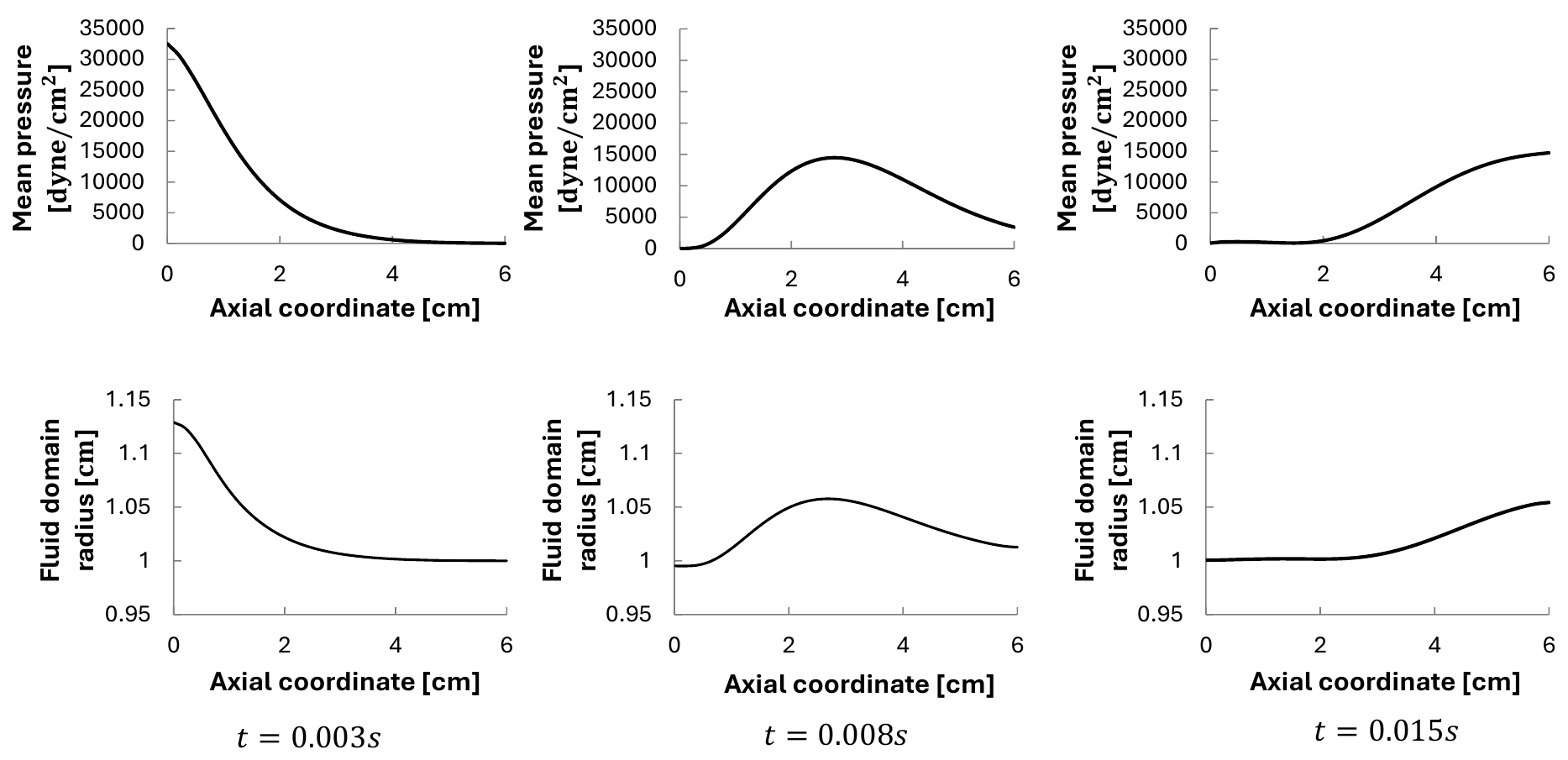}
  \caption{SC-DN-$\alpha$ scheme. Mean pressure (top), obtained as average on the radial section of the fluid domain, and fluid domain radius (bottom) at three instants. Test I.}
  \label{fig:sol_implicito}
\end{figure}

\paragraph{\textbf{Test I}:}
In this first test, we impose the pressure profile
\[
p_{in} = \begin{cases}
    
4\times 10^4\left(1 - \cos\left( \frac{\pi t}{2.5\times10^{-3}}\right)\right) \, dyne/cm^2 \qquad t\le 0.005s, \\
0 \, dyne/cm^2 \qquad\qquad\qquad\qquad\qquad\qquad 0.005s<t<T,
\end{cases}
\]
as a Neumann boundary condition at the inlet section $\Gamma_{\text in}$ and absorbing conditions at the outlet section $\Gamma_{\text out}$ \cite{nobile2008effective,nobile2014inexact}. 
Regarding the choice of the parameter $\alpha$ in this test, we manually tune it to select the optimal one. 
Table \ref{tab:alpha_iterations} shows the average number of iterations for the first 30 time steps employed by the SD-LC-$\alpha$ algorithm, for different values of $\alpha$ ensuring convergence. 
We observe that values of $\alpha$ larger than $0.15$ (thus in particular for standard SC-DN) caused the simulation to diverge. 
This behavior confirms the presence of an upper bound in $\alpha$ to ensure convergence as highlighted by Proposition 1. 
The results of the simulations are presented in Figure \ref{fig:sol_implicito}, for three time frames, showing the axial propagation of the impulse input along the domain. 

\begin{table}[htpb]
\footnotesize
\caption{Average number of iterations over 30 time steps for different values of the parameter $\alpha$. Test I.}\label{tab:alpha_iterations}
\begin{center}
\rowcolors{2}{gray!15}{white}
\begin{tabular}{c|c}
\noalign{\smallskip}\hline\noalign{\smallskip}
    $\alpha$ & iterations  \\
\noalign{\smallskip}\hline\noalign{\smallskip}
    0.05 & 181 \\
    0.10  & 89 \\
    0.13 & 68 \\
    0.14 & 63 \\
    0.15 & 58 \\
    0.16 & X \\
    1.00 (SC-DN) & X\\
    \hline
\end{tabular}
\end{center}  
\end{table}

\paragraph{\textbf{Test II}:}
Given the same input pressure profile as in Test I, in this test we explore the dependence of the optimal value of $\alpha$ on the density ratio, i.e., on the amount of added mass.
In particular, from Table \ref{tab:alpha_density}, we observe an almost constant value of the number of iterations with respect to $\displaystyle{{\rho_s}/{\rho_f}}$. 
In addition, the best value of $\alpha$ seems itself to be quite independent of such a ratio. 
Thus, we conclude that the SC-DN-$\alpha$ is very robust with respect to the density ratio and is able to effectively mitigate the growing convergence difficulties encountered as the added mass increases \cite{causin2005added,forster2007artificial}.

\begin{table}[htpb]
    \footnotesize
    \caption{Dependence of the optimal relaxation parameter $\alpha$ and the number of iterations on the density ratio ${\rho_sh_s}/{\rho_f\mu_i}$. The number of iterations is averaged over 30 time steps. Test II.}\label{tab:alpha_density}
    \begin{center}
    \rowcolors{2}{gray!15}{white}
    \begin{tabular}{l|c|c}
\noalign{\smallskip}\hline\noalign{\smallskip}
         $\displaystyle{\frac{\rho_s}{\rho_f}}$ & $\alpha$ & iterations \\
\noalign{\smallskip}\hline\noalign{\smallskip}
            1.100 & 0.1500 & 58 \\
            0.500 & 0.1450 & 60 \\
            0.100 & 0.1350 & 65 \\
            0.010 & 0.1345 & 65 \\
            0.005 & 0.1335 & 66 \\
         \hline
    \end{tabular}
    \end{center}
\end{table}

\paragraph{\textbf{Test III:}}
In this test, we analyze the convergence behavior of the SC-DN-$\alpha$ algorithm with respect to variations of other physical parameters. 
We impose ramp-shaped pressure at the inlet: 
\[
p_{\text{in, ramp}} =
    \begin{cases}
        10^4\left(1 - \cos\left(\dfrac{\pi t}{2.5\times 10^{-3}}\right)\right)\, \text{dyne}/\text{cm}^2 & t \le 0.0025\,\text{s}, \\
        10^4\, \text{dyne}/\text{cm}^2 & 0.0025\,\text{s} < t \le T.
    \end{cases}
\]
For each scenario, we compare the average number of iterations leading to convergence performed by considering two different selection strategies of $\alpha$: manual selection of the optimal value ensuring convergence (SC-DN-$\alpha_m$) and selection by considering the optimal value estimated in the context of {\sl Minimum-Residual} iterative method (SC-DN-$\alpha_{MR}$) used as a specific Richardson update to obtain our new scheme, see Section \ref{subsec:SC-DN-alpha-scheme}. 
In Table \ref{tab:iterations_ramp} we report the results.

Although the application of the MR update of $\alpha$ to the FSI linear system is not rigorously justified since the monolithic matrix is non-symmetric, these results show its effectiveness in reducing the computational cost in all scenarios with respect to the manual constant selection, thus representing a starting point to devise more rigorous and even efficient update rules for $\alpha$ when dealing with non-symmetric and non-definite-positive matrices.

\begin{table}[htpb]
\footnotesize
\caption{Average number of iterations per time step obtained with the SC-DN-$\alpha$ method for different strategies to select the parameter $\alpha$: manual selection of the optimal value, SC-DN-$\alpha_m$, and automatic selection through the Minumum-Residual update, SC-DN-$\alpha_{MR}$.
For the Lamé constants we consider multiple values of the reference ones $c$ and $\lambda$. REF = reference scenario. Test III.}\label{tab:iterations_ramp}
\begin{center}
\rowcolors{2}{gray!15}{white} 
\begin{tabular}{l|c|c}
\noalign{\smallskip}\hline\noalign{\smallskip}
Scenarios & SC-DN-$\alpha_{m}$ & SC-DN-$\alpha_{MR}$ \\
\noalign{\smallskip}\hline\noalign{\smallskip}

REF                & 47.9 & 22.2 \\
$\rho_s = 0.01$   & 42.2 & 24.6 \\
$\Delta t = 5e\!-\!4$   & 87.9 & 39.1 \\
$5c-5\lambda$     & 53.3 & 21.2 \\
$c/5-\lambda/5$    & 47.8 & 22.7 \\
$H_s = 0.15$      & 42.1 & 18.0 \\
$H_s = 0.05$      & 63.6 & 35.8 \\
\hline
\end{tabular}
\end{center}
\end{table}

\section{A new loosely coupled partitioned scheme: LC-DN-$\alpha$ algorithm}\label{sec:DN-alpha-LC}
Given the promising results obtained for the SC-DN-$\alpha$ scheme, in this Section we present a new loosely-coupled scheme derived as a single preconditioned Richardson iteration $P \mathbf{X}^{n+1} = P \mathbf{X}^{n} + \alpha \left(\mathbf{b}^{n+1} - A \mathbf{X}^{n}  \right)$,  with the same preconditioner $P$ as the one described in Section \ref{subsec:SC-DN-alpha-scheme}.
However, it can be easily verified that the resulting scheme is not consistent. 
In particular, we obtain the following asymptotic truncation errors: 
\begin{align*}
    &\boldsymbol{\tau}^{n+1}_f = - \rho_f(1 -\alpha) \frac{\bm{u}(t^n) - \bm{u}(t^{n-1})}{\Delta t}, \\
    &\boldsymbol{\tau}^{n+1}_s = -\rho_sh_s (1-\alpha)\left(\frac{2\boldsymbol{\eta}(t^n) -3\bm{\eta}(t^{n-1}) +\bm{\eta}(t^{n-2})}{\Delta t^2}\right),\\
    & \boldsymbol{\tau}^{n+1}_{\bm w} = (1 -\alpha)\bm{w}(t^n)
\end{align*}
where $\bm{\tau}^{n+1}_f$ and $\bm{\tau}^{n+1}_s$ are associated with the fluid and structure momentum problems, and $\bm{\tau}^{n+1}_{\bm w}$ is associated with the structure velocity update.

To restore consistency in the LC scheme, we explicitly subtract the discrete counterparts of these errors from the governing equations, obtaining a novel algorithm, referred to as LC-DN-$\alpha$ (see Algorithm \ref{alg:DN-alpha-LC}).

\begin{algorithm}
\caption{LC-DN-$\alpha$ scheme}
\label{alg:DN-alpha-LC}
{\sl Given suitable initial conditions and boundary conditions on $\partial\Omega_f^t\setminus\Sigma^t$ and $\partial\Omega_s^t\setminus\Sigma^t$, find $\forall n = 0,1, \dots, T/ \Delta t$, $\bm{u}^{n+1}, p^{n+1}$, ${\boldsymbol{\eta}}^{n+1}$, and $\bm{w}^{n+1}$ by solving in sequence ($n+1$ understood):}
\vspace{1pt}
\begin{enumerate}[leftmargin=*, align=left]
    \item Fluid problem on the unknowns $\bm{u}^{n+1} \text{ and } p^{n+1}$ equipped with a Dirichlet boundary condition at the interface $\Sigma$: 
    \begin{align}
    \left\{
        \begin{aligned}
            &\rho_f \left(\frac{\bm{u}}{\Delta t} + (\bm{u}^n - \bm{u}_{ALE}^n) \cdot \nabla \bm{u}\right)- \nabla \cdot \boldsymbol{T}_{f} =  \rho_f \frac{\bm{u}^n}{\Delta t} + \notag \\ 
		& \qquad \qquad + (1- \alpha) \left[\rho_f \left(\frac{\bm{u}^{n} -   \bm{u}^{n-1}}{\Delta t} + (\bm{u}^{n} - \bm{u}_{ALE}^{n}) \cdot \nabla \bm{u}^{n}\right) - \nabla \cdot \boldsymbol{T}_{f}^{n}\right]
            &{} \text{in}\,\, \Omega_f,\\[6pt]
            &\nabla \cdot \bm{u} = (1 - \alpha)\nabla \cdot \bm{u}^{n}
            &{} \text{in}\,\, \Omega_f,\\[6pt]
                &\bm{u} = (1 - \alpha)\bm{u}^{n} + \alpha \bm{w}^n
            &{} \text{on}\,\, \Sigma;
        \end{aligned}
    \right.
\end{align}

    \item Structure problem on the unknown $\boldsymbol{\eta}^{n+1}$ equipped with a Neumann boundary condition at the interface $\Sigma$: 
    \begin{align}
    \left\{
    \begin{aligned}
        &\rho_s \frac{{\boldsymbol{\eta}}}{\Delta t^2} - \nabla \cdot {\boldsymbol{T}}_{s}  + \beta\bm{\eta}=  \rho_s\frac{2{\boldsymbol{\eta}}^n - {\boldsymbol{\eta}}^{n-1}}{\Delta t^2} + \\
        & \qquad \qquad \qquad+ (1-\alpha) \left( \rho_s \frac{{\boldsymbol{\eta}}^{n} - 2{\boldsymbol{\eta}}^{n-1} + {\boldsymbol{\eta}}^{n-2} }{\Delta t^2} - \nabla \cdot {\boldsymbol{T}}_{s}^n\right)
        &{} \text{in } {\Omega}_s^0,\notag \\[6pt]
        &\boldsymbol{T}_{s} \cdot \bm{n} = -  \boldsymbol{T}_{f} \cdot \bm{n} + (1-\alpha) \left(\boldsymbol{{T}}_{s}^n \cdot \bm{n}+ \boldsymbol{T}_{f}^n \cdot \bm{n} \right)
        &{} \text{on}\ \Sigma^0; \notag 
    \end{aligned}
    \right.
    \end{align}
    \item Update $\displaystyle{\bm{w} = \frac{\bm{\eta} - \bm{\eta}^n}{\Delta t}\bigg|_{\Sigma^0}}$.
    \end{enumerate}
\end{algorithm}

\subsection{Stability of the LC-DN-$\alpha$ algorithm}\label{subsec:stability}

For the study of the stability properties of Algorithm \ref{alg:DN-alpha-LC}, we rely on the same benchmark problem described in Section \ref{subsec:convergence}. 
By considering again the backwards Euler and the implicit BDF1 schemes for the approximation in time, Algorithm \ref{alg:DN-alpha-LC} applied to the resulting discretized-in-time benchmark problem states as follows: 

{\sl Find $\bm{u}^{n+1}$, $p^{n+1}$, ${\eta}^{n+1}$, and $w^{n+1}$ by solving in sequence ($n+1$ understood):}
\begin{itemize}
    \item[(a)] Fluid problem equipped with Dirichlet boundary condition at the interface $\Sigma$:
        \begin{align}\label{eq:model_LC_fluid}
            \left\{
            \begin{aligned}
                &\rho_f \frac{\bm{u}}{\Delta t} + \nabla p = \rho_f \frac{\bm{u}^{n}}{\Delta t} +(1-\alpha)\left(\rho_f \frac{\bm{u}^{n}-\bm{u}^{n-1}}{\Delta t} + \nabla p^{n}\right)
                &{} \text{in}\ \Omega_f, \\
	           &\nabla \cdot \bm{u}= (1-\alpha) \nabla \cdot \bm{u}^{n}
                &{} \text{in}\ \Omega_f, \\
	           &\bm{u} \cdot \bm{n} = (1-\alpha) \bm{u}^{n} \cdot \bm{n} + \alpha w^n
                &{} \text{on}\ \Sigma;
            \end{aligned}
            \right.
        \end{align}
    \item[(b)] Structure problem:
        \begin{align}\label{eq:model_LC_structure}
            & \rho_s h_s \frac{\eta}{\Delta t^2} + \beta \eta -  p = \rho_s h_s\frac{2\eta^{n} - \eta^{n-1}}{\Delta t^2} + \notag\\ 
            & \qquad + (1-\alpha)\left(\rho_s h_s \frac{\eta^{n}-2\eta^{n-1} + \eta^{n-2}}{\Delta t^2} + \beta \eta^n -  p^{n}\right) 
            &{} \text{in } \Omega_s.
        \end{align} 
    \item[(c)] Update $\displaystyle{{w} = \frac{\eta - \eta^n}{\Delta t }}$.
\end{itemize}

In the following two propositions, we provide sufficient conditions which ensure absolute instability and stability, respectively, of the LC-DN-$\alpha$ scheme in the range $\alpha \in (0,1]$.
 
\begin{proposition}\label{prop:LC-unstability}
Assume $0<\alpha\leq1$; then the LC-DN-$\alpha$ scheme is unconditionally unstable, when 
\begin{equation}\label{prpo:hyp}
    \rho_f\mu_i > \frac{\beta(2-\alpha)\Delta t^2}{4\alpha} + \left(\frac{2(1-\alpha)}{\alpha} + 1\right)\rho_sh_s \quad \forall i .
\end{equation} 
\end{proposition}

\begin{proof}
The same reasoning as in Section \ref{subsec:convergence} applied to \eqref{eq:model_LC_fluid}-\eqref{eq:model_LC_structure} yields to the following finite difference equation:
\begin{align*}
    &\rho_s h_s \frac{\eta_i^{n+1} - 2 \eta_i^n + \eta_i^{n-1}}{\Delta t^2} 
    + \beta \left(\eta_i^{n+1}-(1-\alpha)\eta_i^n\right) +\notag \\
    & \qquad\qquad - (1-\alpha ) \left( \rho_s h_s \frac{\eta_i^{n} - 2 
    \eta_i^{n-1} + \eta_i^{n-2}}{\Delta t^2} \right) + \notag\\
     & \qquad\qquad + \alpha  \rho_f \mu_i\left(\frac{\eta_i^{n} - 2 \eta_i^{n-1} + \eta_i^{n-2}}{\Delta t^2} \right) = 0
    \,\, \text{on } \Sigma.
\end{align*}
We introduce its associated characteristic polynomial $\Chi(y)\in \mathbb{P}^3$:

\begin{align}\label{eq:polychar_B2}
    & \Chi(y) = \left[\frac{\rho_sh_s}{\Delta t^2} + \beta\right]y^3 
    + \left[\alpha\frac{\rho_f\mu_i}{\Delta t^2} -\beta(1-\alpha) - 3\frac{\rho_sh_s}{\Delta t^2} +\alpha\frac{\rho_sh_s}{\Delta t^2}\right]y^2 +\notag\\
    & \qquad\qquad + \left[3\frac{\rho_sh_s}{\Delta t^2} -2\alpha\frac{\rho_f\mu_i}{\Delta t^2} -2\alpha\frac{\rho_sh_s}{\Delta t^2}\right]y + \left[-(1-\alpha)\frac{\rho_sh_s}{\Delta t^2} + \alpha\frac{\rho_f\mu_i}{\Delta t^2}\right].
\end{align}
Then, we evaluate this polynomial at ${y}=-1$: 
\begin{equation*}
    \Chi(-1) = \frac{\beta\Delta t^2(\alpha-2) +4\alpha(\rho_sh_s + \rho_f\mu_i) - 8\rho_sh_s}{\Delta t ^2},
\end{equation*}
which is positive under the hypothesis \eqref{prpo:hyp}.
Since $\displaystyle{\lim_{y\rightarrow -\infty}{\Chi(y)}} = -\infty$ and $\displaystyle{\lim_{y\rightarrow +\infty}{\Chi(y)}} = +\infty$, it follows that under such hypothesis the characteristic polynomial has at least one real root with modulus greater than 1. 
This implies that the scheme is absolutely unstable. $\qquad \square$

\end{proof}

\begin{remark}
    In view of \eqref{prpo:hyp}, it appears that decreasing $\alpha$ amplifies the right-hand side contribution, effectively shifting the instability threshold for $\rho_f$ to higher values; this makes it more difficult for $\rho_f$ to exceed such a threshold.
    We conclude that decreasing the value of $\alpha \in (0,1)$ allows us to decrease the range of unstable $\rho_f$ and thus, in principle, to consider larger added mass than for $\alpha=1$. 
\end{remark}

\begin{remark}
    From \eqref{prpo:hyp}, observe that for $\alpha=1$ we recover the same instability condition of \cite{causin2005added}, adapted to the case of BDF1 instead of Leap-Frog for the time discretization of the structure problem.
\end{remark}

\begin{proposition}\label{prop:3}
Assume $0<\alpha\leq1$ and $\frac{\rho_s h_s}{\rho_f\mu_i} < 1$; then the LC-DN-$\alpha$ scheme is unconditionally absolute stable if, for any $i$,
\begin{equation}\label{prop:unconditional_stability}
    0 < \alpha < \overline{\alpha} = \frac{2\frac{\rho_sh_s}{\rho_f\mu_i} }{\frac{\rho_sh_s}{\rho_f\mu_i} + 1} <1.
\end{equation}

\end{proposition}

\begin{proof}
We start by introducing the polynomial ${\mathcal{Q}(y):= \frac{\Delta t^2}{\rho_sh_s}\mathcal{\Chi}(y)} = a_3 y^3 + a_2 y^2 + a_1 y + a_0$ , where $\Chi(y)$ is the characteristic polynomial \eqref{eq:polychar_B2} and the coefficients $a_{i}$ are defined as: 
\begin{align*}
    a_0 &= \alpha\frac{\rho_f\mu_i}{\rho_sh_s} - 1 +\alpha, \\
    a_1 & = 2  -2\alpha -2\alpha\frac{\rho_f\mu_i}{\rho_sh_s},\\
    a_2 & = -3 +\alpha -(1-\alpha)\frac{\beta\Delta t^2}{\rho_sh_s} +\alpha\frac{\rho_f\mu_i}{\rho_sh_s},\\
    a_3 & = 1 + \frac{\beta\Delta t^2}{\rho_sh_s}.
\end{align*}

To analyze the stability of the scheme \eqref{eq:model_LC_fluid}-\eqref{eq:model_LC_structure}, we verify the necessary and sufficient conditions given by the Jury's criterion \cite{jury1963roots} in the regime ${\frac{\rho_sh_s}{\rho_f\mu_i}<1}$: 
\begin{enumerate}
    \item $\mathcal{Q}(1)>0$;
    \item $\mathcal{Q}(-1)<0$;
    \item $|a_3|>|a_0|$;
    \item $|a_0^2 -a_3^2|> |a_0a_2-a_1a_3|$.
\end{enumerate}

Evaluating $\mathcal{Q}(y)$ at ${y}=1$, readily shows that Condition 1 is satisfied for any value of $\Delta t$ provided that $\alpha>0$. 

Next, consider $\mathcal{Q}(-1)$: 
\begin{align*}
    \mathcal{Q}(-1) = \frac{-8\rho_sh_s + 4\alpha(\rho_sh_s +\rho_f\mu_i) +\beta\Delta t^2(\alpha-2)}{\Delta t^2}.
\end{align*}
A straightforward calculation shows that, for $\alpha<\overline{\alpha}$, the latter is negative for any $i$, and thus  Condition 2 is fulfilled for any $\Delta t$.

Condition 3 reads as follows: 
\begin{equation*}
     \left|1 + \frac{\beta\Delta t^2}{\rho_sh_s}\right| > \left|\alpha\frac{\rho_f\mu_i}{\rho_sh_s} - 1 +\alpha\right|,
\end{equation*}
which is easy to show that it is satisfied for any $i$ and $\Delta t$ for example for $\alpha<\overline{\alpha}$.

Finally, to assess Condition 4, note that it can be equivalently written as
\begin{align} \label{eq:zeta-psi-relations}
    |\zeta_i| > |\zeta_i +\psi_i| &\iff 
    \begin{cases}
        \zeta_i> 0 \\
        \psi_i<0 
    \end{cases}
    \bigcup\,\,\,\,
    \begin{cases}
        \zeta_i < 0 \\
        \psi_i>0
    \end{cases} \quad \forall i,
\end{align}
where
\begin{align*}
    \zeta_i(\alpha) & = \left(\alpha\frac{\rho_f\mu_i}{\rho_sh_s} - 1 + \alpha\right)^2 - \left(1 + \frac{\beta\Delta t^2}{\rho_sh_s}\right)^2, \\
    \psi_i(\alpha) & = \alpha^2\left(\frac{\beta\Delta t^2\rho_f\mu_i}{(\rho_sh_s)^2} + \frac{\beta\Delta t^2}{\rho_sh_s}\right) +\alpha\frac{\beta\Delta t^2\rho_f\mu_i}{(\rho_sh_s)^2} + \left(\frac{\beta\Delta t^2}{\rho_sh_s}\right)^2.
\end{align*}
The roots of $\zeta_i(\alpha)$ can be easily computed and show that $\zeta_i(\alpha)>0$ for $\alpha\ge \alpha_i = \frac{2\rho_sh_s + \beta\Delta t^2}{\rho_sh_s + \rho_f\mu_i}$, while $\psi_i(\alpha)> 0$ for all $\alpha>0$. 
Therefore, \eqref{eq:zeta-psi-relations} is satisfied for ${0<\alpha\leq\alpha_i}$. 
Observing that $\alpha_i>\overline{\alpha}$ for any $i$ and combining all the conditions found from the beginning of the proof, we found that if hypothesis \eqref{prop:unconditional_stability} holds true, then the Jury's criterion is satisfied for any $\Delta t$ and the thesis follows. $\qquad \square$
\end{proof}

\begin{remark}
    Note that condition \eqref{prop:unconditional_stability} can be formulated also as $\rho_f\mu_i< \frac{2-\alpha}{\alpha}\rho_sh_s$. We observe that the right-hand side of this inequality increases for decreasing $\alpha$. 
    Thus, small values of $\alpha$ allow us to extend the range of added mass for which the LC-DN-$\alpha$ scheme is stable, in principle covering all possible added mass regimes provided that $\alpha$ is small enough.
    
\end{remark}

\subsection{Numerical results}\label{subsec:resultsLC}
We now present a set of 2D numerical experiments performed with the LC-DN-$\alpha$ scheme (Algorithm \ref{alg:DN-alpha-LC}) with realistic equations, in order to validate the theoretical results introduced in the previous section.
These tests are carried out considering the problem presented in Section \ref{subsec:convergence-results}, using the same discretization setting and physical parameters, unless otherwise specified. 
In all scenarios, we impose the following inlet pressure profile:
\[
p_{in} = \begin{cases}
    
2\times 10^4\left(1 - \cos\left( \frac{\pi t}{0.05}\right)\right) \, dyne/cm^2 \qquad t\le 0.1s, \\
0 \, dyne/cm^2 \qquad\qquad\qquad\qquad\qquad\qquad 0.1s<t<T.
\end{cases}
\]

For the sake of comparison, alongside the LC-DN-$\alpha$ solution, we also report the SC-DN-$\alpha$ one, with $\alpha$ chosen from the values that guarantee convergence.
It should be noted that the convergent SC solution is independent of $\alpha$.

\begin{figure}[ht]
  \centering
  \includegraphics[width=1\linewidth]{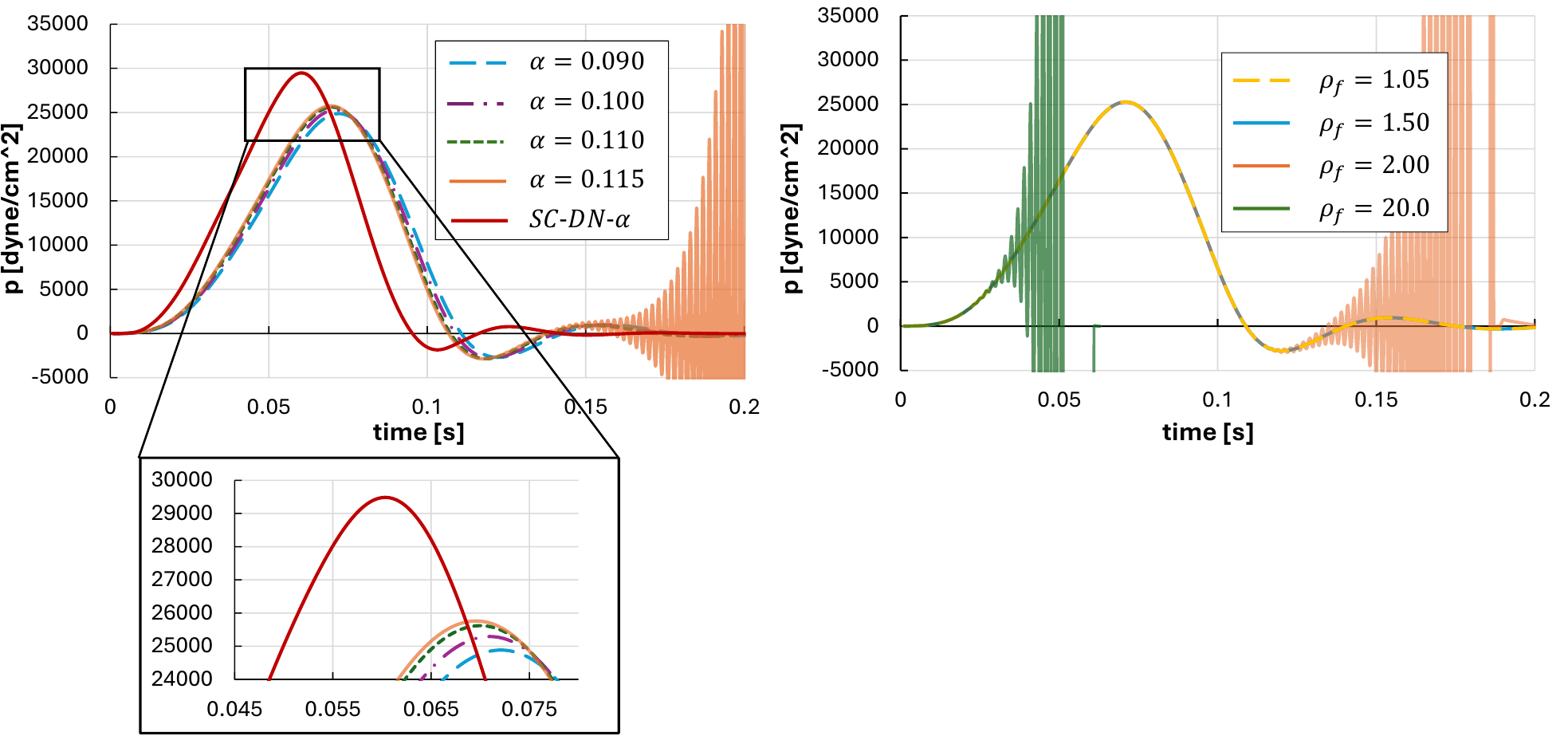}
  \caption{Fluid mean pressure over the cross-section at $z = 3\,cm$: different values of $\alpha$ for $\rho_s/\rho_f = 0.9$ (left); different values of $\rho_f$ for $\rho_s = 1$ and $\alpha=0.11$ (right). $\Delta t=10^{-3}s$. Test IV.}
  \label{fig:conditional-stabilityalpha}
\end{figure}

\paragraph{\textbf{Test IV:}}
In this test, first we perform numerical experiments with different values of $\alpha$ and fixed $\Delta t=1.0\times10^{-3}s$ and $\rho_s/\rho_f = 0.9$ (thus in the presence of significant added mass).
The results are presented in Figure \ref{fig:conditional-stabilityalpha} - left.
We observe that if $\alpha$ is small enough, then the LC-DN-$\alpha$ solutions are stable and slightly delayed and smoothed with respect to the SC-DN-$\alpha$ one. 
We also notice that, for such values of $\alpha$, the LC-DN-$\alpha$ solution approaches the SC-DN-$\alpha$ one for increasing values of $\alpha$.
Moreover, we notice that if $\alpha$ is greater than a certain threshold, the solution of the LC-DN-$\alpha$ becomes unstable (orange plot in Figure \ref{fig:conditional-stabilityalpha}). 
This confirms the theoretical results of Prepositions 2 and 3, as higher values of this parameter shrink the stability domain.

Secondly, we perform numerical experiments with different values of $\rho_f$ and fixed $\Delta t = 10^{-3}s$, $\alpha= 0.11$, and $\rho_s=1\,g/cm^3$ (see Figure \ref{fig:conditional-stabilityalpha} - right).
We notice that for the two smallest values of $\rho_f$ the solution is stable, for increasing values of $\rho_f$ the solution becomes unstable. 
This agrees with the findings of Proposition 2, which identifies a threshold for $\rho_f$ over which the solution is expected to be unconditionally stable.
Interestingly, for $\rho_f=1.5g/cm^3$, $\Delta t = 10^{-4}s$, $\alpha= 0.11$, and $\rho_s=1\,g/cm^3$ we found an unstable solution (not reported here for the sake of exposition), confirming that the threshold of Proposition 2 decreases for decreasing $\Delta t$ (notice that for the same parameters and $\Delta t=10^{-3}s$ the solution was stable).

\begin{figure}[ht]
  \centering
  \includegraphics[width=1\linewidth]{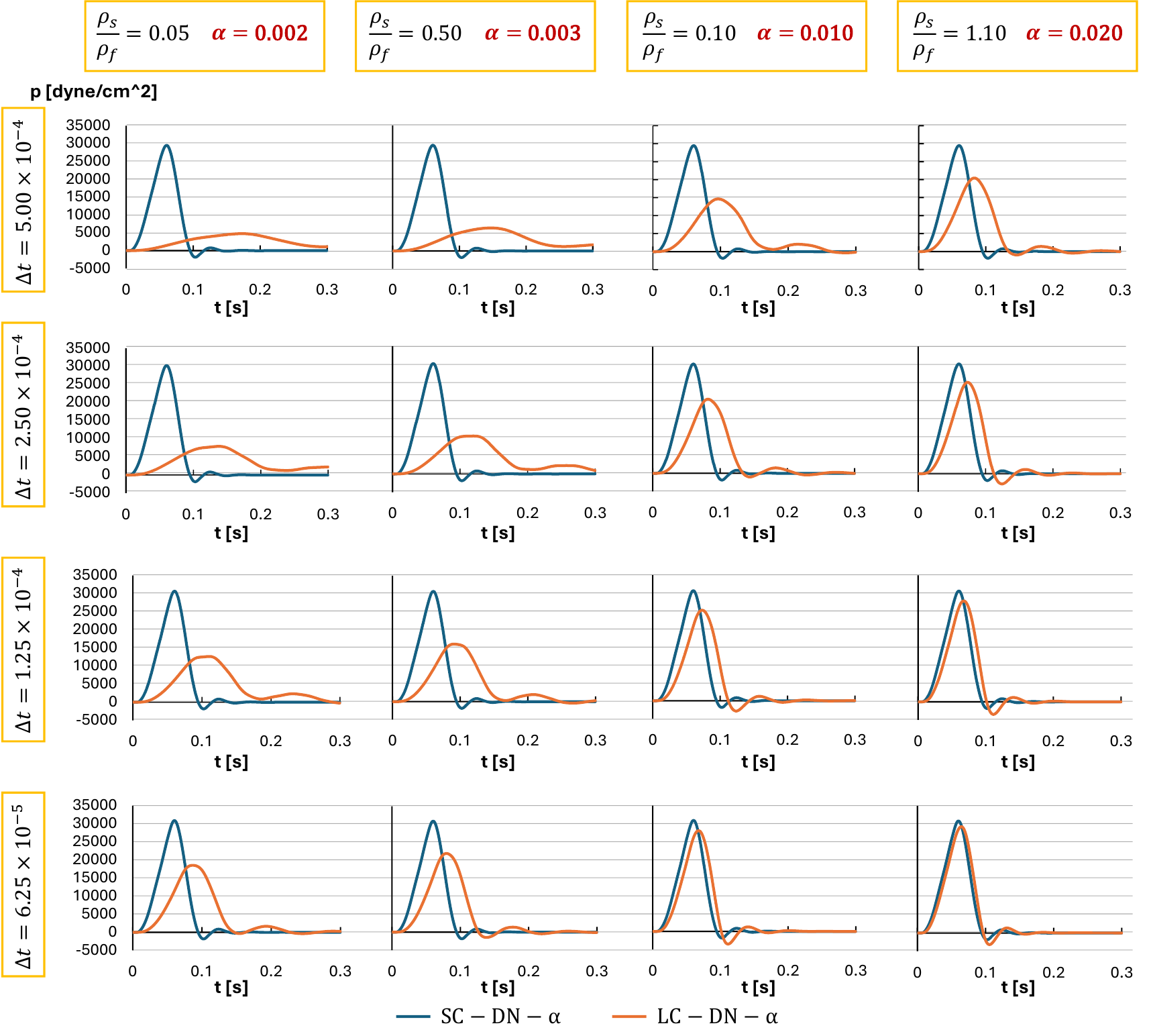}
  \caption{Comparison between the SC-DN-$\alpha$ and the LC-DN-$\alpha$ fluid mean pressure over the cross-section at $z = 3\,cm$, for different values of $\rho_s/\rho_f$ and $\Delta t$. Test V.}
  \label{fig:convergence}
\end{figure}

\paragraph{\textbf{Test V:}}
In the last test, we perform simulations for decreasing density ratios and time steps, see Figure \ref{fig:convergence}.

Each column of the Figure shows that the LC-DN-$\alpha$ solution gradually converges towards the SC-DN-$\alpha$ for decreasing values of $\Delta t$, demonstrating that the LC-DN-$\alpha$ scheme is consistent.

Secondly, each row of the Figure shows that one needs smaller values of $\alpha$ to be stable when $\rho_s/\rho_f$ decreases, in accordance with the value of $\overline{\alpha}$ in Proposition 3.
Each row also shows that accuracy improves as the added mass decreases, provided that the value of $\alpha$ is large enough.

\vspace{5pt}
In conclusion, these findings of Tests IV and V confirm the consistency and stability when $\alpha$ is suitably chosen of the LC-DN-$\alpha$ scheme with respect to the added-mass effect, highlighting its potential for applications, e.g., in hemodynamics.
Naturally, its performance in more realistic scenarios  remains to be explored and will be the subject of future work.

\clearpage

\section*{Acknowledgements}
FR, CV acknowledge their membership in INdAM GNCS - Gruppo Nazionale per
il Calcolo Scientifico (National Group for Scientific Computing, Italy). 
FR has been partially supported by the INdAM GNCS project 2025 - CUP E53C24001950001.
FR has been funded by the European Union-Next Generation EU, Mission 4, Component 1, CUP: D53D23018770001, under the research project MIUR PRIN22-PNRR n.P20223KSS2, \lq\lq Machine learning for fluid structure interaction in cardiovascular problems: efficient solutions, model reduction, inverse problems\rq\rq.
CV has been partially supported by i) the European Union-Next Generation EU, Mission 4, Component 1, CUP: D53D23018770001, under the research project MIUR PRIN22-PNRR n.P20223KSS2, \lq\lq Machine learning for fluid structure interaction in cardiovascular problems: efficient solutions, model reduction, inverse problems\rq\rq; ii) the Italian Ministry of Health within the PNC PROGETTO HUB LIFE SCIENCE - DIAGNOSTICA AVANZATA (HLSDA) \lq\lq INNOVA\rq\rq, PNCE3 - 2022-23683266-CUP: D43C22004930001, within the \lq\lq Piano Nazionale Complementare Ecosistema Innovativo della Salute\rq\rq - Codice univoco investimento: PNCE3- 2022-23683266; iii) the Italian research project MIUR PRIN22 n.2022L3JC5T \lq\lq Predicting the outcome of endovascular repair for thoracic aortic aneurysms: analysis of fluid dynamic modeling in different anatomical settings and clinical validation\rq\rq; iv) Italian Ministry of Health within the project \lq\lq CAL.HUB.RIA\rq\rq - CALABRIA HUB PER RICERCA INNOVATIVA ED AVANZATA. Code: T4-AN-09, \\ CUP: F63C22000530001.


%
%

\bibliographystyle{spmpsci}      
\bibliography{./Renzi-et-al-JSC-bibliography}   

\end{document}